\documentclass[a4paper,12pt]{article}

\usepackage{amsmath, amssymb, amsthm}
\usepackage{mathtools}
\usepackage[dvips,final]{graphicx}

\usepackage{hyperref}
\hypersetup{
	colorlinks=true,
	linkcolor=blue,
	urlcolor=cyan,
	citecolor=blue,
}

\usepackage{color}

\definecolor{gre}{rgb}{.2,.7,0}
\definecolor{ora}{rgb}{1,0.6,0}
\definecolor{blue1}{rgb}{0,.3,.6}
\definecolor{gbl}{rgb}{0,.4,.4}
\definecolor{shi}{rgb}{.5,.25,0}

\usepackage{titlesec}
\titleformat{\section}{\bfseries}{\thesection}{1em}{}
\titleformat{\subsection}{\itshape}{\thesubsection}{1em}{}

\numberwithin{equation}{section}

\newfont{\ctv}{msam10}

\newcommand{\bbox}{\mbox{\ctv \symbol{4}}}
\def\QED{{${}\hfill\bbox$}}
\newenvironment{pf}[1]{\par\vskip1mm{\noindent\it #1.}\ }{\QED\par
\vskip2mm}

\def\bpf{\begin{pf}}
\def\epf{\end{pf}}

\def\expe{\hbox{\rm e}}
\def\ve{\varepsilon}

\def\vrt{\phi}

\def\dd{\,\mathrm{d}}
\def\dive{\mathrm{\,div\,}}
\def\sign{\mathrm{\,sign}}

\def\supess{\mathop{\mathrm{\,sup\,ess}}}
\def\for{\mathrm{\ for\ }}
\def\ale{\mathrm{\ a.\,e.}}

\def\play{\mathfrak{p}}
\def\uu{U}
\def\bsi{\boldsymbol{e}}
\def\usi{u_{\bsi}}
\def\bsa{\boldsymbol{a}}
\def\bsb{\boldsymbol{b}}

\def\sumiz{\sum_{i=0}^n}
\def\sumim{\sum_{i=0}^{n-1}}

\def\on{^{(n)}}

\def\real{\mathbb{R}}
\def\nat{\mathbb{N}}

\def\io{\int_{\Omega}}
\def\ipo{\int_{\partial\Omega}}

\def\irn{\int_{\real^N}}

\def\be{\begin{equation}\label}
\def\ee{\end{equation}}

\def\ber{\begin{eqnarray}}
\def\eer{\end{eqnarray}}

\def\bers{\begin{eqnarray*}}
\def\eers{\end{eqnarray*}}

\def\bpf{\begin{pf}}
\def\epf{\end{pf}}

\newtheorem{theorem}{Theorem}[section]
\newtheorem{lemma}[theorem]{Lemma}

\newtheorem{hypothesis}[theorem]{Hypothesis}
\newtheorem{proposition}[theorem]{Proposition}
\newtheorem{definition}[theorem]{Definition}

\begin{document}

\title{Diffusion in porous media with hysteresis\\and bounded speed of propagation
\thanks{The support from the European Union's Horizon Europe research and innovation programme under the Marie Sk\l odowska-Curie grant agreement No 101102708, from the M\v{S}MT grant 8X23001, and from the GA\v CR project 24-10586S is gratefully acknowledged.}
}

\author{Chiara Gavioli
\thanks{Faculty of Civil Engineering, Czech Technical University, Th\'akurova 7, CZ-16629 Praha 6, Czech Republic, E-mail: {\tt chiara.gavioli@cvut.cz}.}
\and Pavel Krej\v c\'{\i}
\thanks{Faculty of Civil Engineering, Czech Technical University, Th\'akurova 7, CZ-16629 Praha 6, Czech Republic, E-mail: {\tt Pavel.Krejci@cvut.cz}.}
\thanks{Institute of Mathematics, Czech Academy of Sciences, {\v{Z}}itn{\'a} 25, CZ-11567 Praha 1, Czech Republic, E-mail: {\tt krejci@math.cas.cz}.}
}

\date{}

\maketitle

\begin{abstract}
It is shown that the problem of moisture propagation in porous media with a nonlinear relation between the mass flux and the pressure gradient as a counterpart of the Darcy law exhibits the property of bounded speed of propagation even in the case of a hysteresis relation between the capillary pressure and the moisture content. The paper specifies conditions for existence and uniqueness of solutions, and provides an upper bound for the moisture propagation speed.

\bigskip


\medskip

\end{abstract}


\section*{Introduction}\label{sec:intr}

\begin{figure}[htb]
	\begin{center}
		\includegraphics[width=9cm]{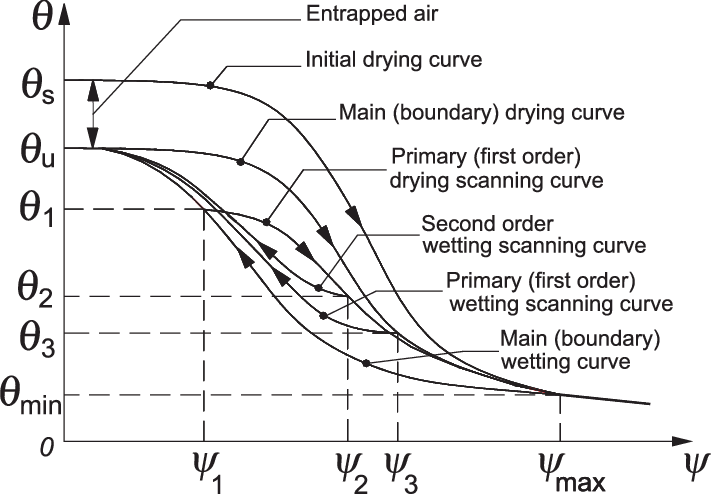}
		\caption{Typical experimental hysteresis dependence in porous media between the logarithm soil suction $\psi$ (which can be interpreted as a decreasing function of the pressure) and the volumetric water content $\theta$.}
		\label{f1}
	\end{center}
\end{figure}

The problem of bounded speed of propagation in classical problems of diffusion with no mass exchange on the boundary of the spatial domain has been extensively studied in the mathematical literature, see \cite{vazquez} and the references therein. To our knowledge, nothing has been done so far in this direction if hysteresis effects in the constitutive law are taken into account. Note that in porous media, hysteresis occurs as a result of surface tension on the liquid-gas interface, and in many cases cannot be ignored, see \cite{bear}. This paper fills this gap and investigates the problem of bounded propagation speed in presence of hysteresis.
To be precise, we study the following model problem in dimensionless form
\begin{align}\label{e1}
\theta_t &= \dive \big(\kappa(x,\theta) |\nabla u|^{p-2}\nabla u\big),\\ \label{e1a}
\theta &= G[u],
\end{align}
with unknown function $u$ which represents the normalized pressure, that is, $u = (\pi - \pi_0)/\pi_0$, where $\pi$ is the physical pressure and $\pi_0$ is the standard pressure. By $\theta$ we denote the saturation which, up to a linear transformation, can be identified with the moisture content. It is related to the pressure in \eqref{e1a} in terms of a Preisach hysteresis operator $G$ with initial memory $\lambda$. Experimental validation of the Preisach operator for modeling the hysteretic pressure-saturation relationship in porous media was provided in \cite{flynn1,flynn2} through a rigorous fitting procedure. While some models (see e.\,g.\ \cite{beha,gosch}) employ a play operator for hysteresis, this approach is limited to horizontal scanning curves, unlike the behavior shown in Figure~\ref{f1} taken from \cite{pfb}. The Preisach construction, based on averaging a continuum of play elements (see Definition~\ref{dpr} below), more accurately captures the full range of wetting and drying curves. In addition, models using the play operator usually incorporate a rate-dependent regularization in the pressure-saturation relationship. Here we consider instead a purely hysteretic Preisach-type relation. This introduces a significant mathematical challenge: the problem becomes degenerate in the sense that knowing $G[u]_t$ does not give a complete information about $u_t$. More specifically, at every point $x \in \Omega$ and every time $t_0$ where $u_t$ changes sign (a so-called {\em turning point\/}), we have
\be{e6a}
u_t(x,t_0-\delta)\cdot u_t(x,t_0+\delta)<0 \ \ \forall \delta\in (0,\delta_0(x)) \ \Longrightarrow \ \liminf_{\delta \to 0+} \frac{G[u]_t(x,t_0+\delta)}{u_t(x,t_0+\delta)} = 0,
\ee
that is, the starting slope is horizontal. This again corresponds to the behavior shown in Figure~\ref{f1}. We refer to \cite{colli} for more details.

Equation~\eqref{e1}, where $\kappa(x,\theta)$ denotes a saturation-dependent specific permeability, follows from the mass balance and a saturation-dependent variant of the nonlinear Darcy law with exponent $p>2$. In the physics and engineering literature on porous media, the nonlinear Darcy law with mass flux proportional to a power $p>2$ of the pressure gradient is used to model flows at very low velocities (also called \textit{pre-Darcy flow}), see \cite{dejam}. From an analytical point of view, 
in the case of non-degenerate pressure-saturation law, the assumption $p>2$ ensures that the moisture front propagation speed is finite, as shown in \cite[Chapter~VI]{dib}. The consequence is that solutions with compactly supported initial data are compactly supported for any positive time, with a support that expands with time. In problems with doubly degenerate diffusion (or degenerate-singular according to the terminology of \cite{iva,boge}), a bigger exponent $p$ may be necessary, and we postpone the discussion on this issue to Subsection~\ref{nece}. Let us just mention that for Problem~\eqref{e1}--\eqref{e1a}, the existence of a clear interface separating the wet and dry regions, as observed in reality, is guaranteed for $p>3$.

Problem~\eqref{e1}--\eqref{e1a} is considered in a bounded Lipschitzian domain $\Omega \subset\real^N$ and time interval $(0,T)$, and with a given initial condition
\be{e3}
u(x, 0) = u_0(x)
\ee
and a general boundary condition including the homogeneous Neumann, Dirichlet, or Robin conditions, see \eqref{e2}. From the physical point of view, Dirichlet boundary condition is difficult to justify, but it constitutes an intermediate step in the problem of bounded speed of propagation. Let $B_R$ denote the open ball of radius $R$ centered at the origin, and assume that there exist $0<R_0 < R_1$ such that $\overline{B_{R_1}} \subset \Omega$, and the initial condition $u_0$ as well as the initial memory $\lambda$ of $G$ vanish outside $\overline{B_{R_0}} \subset B_{R_1}$, see Definition~\ref{dpr} and Hypothesis~\ref{hR} below. We show in Theorem~\ref{t2} that there exists a time $t_1>0$ such that the solution to \eqref{e1}--\eqref{e3} with homogeneous Dirichlet/Neumann/Robin boundary conditions on $\partial\Omega$ vanishes outside $B_{R_1}$ for $t \in [0,t_1]$. The solution is then extended to a solution on the whole space $\real^N \times (0,\infty)$ with support in a ball $B_{R(t)} \subset \real^N$. An explicit formula for $R(t)$ represents an upper bound on the moisture propagation speed. The proof involves a comparison with a traveling wave solution which can only be constructed if $p>3$. A similar technique has been used also in the case without hysteresis, see \cite{diaz}; see also \cite[Chapter~VI]{dib} for a proof by comparison using the so-called Barenblatt solution which, however, is not suitable for the hysteretic case.

The degeneracy associated with the Preisach operator, which is not considered in the classical theory of parabolic equations with hysteresis developed by Visintin in \cite{vis}, represents a major difficulty in proving the existence of solutions. In order to control the time derivative of $u$, we have to restrict our considerations to the so-called convexifiable Preisach operators following the technique developed in \cite{colli} and extended in \cite{perme} to the case of saturation-dependent permeability. As noted in \cite{colli}, Preisach operators used in engineering practice are typically convexifiable, so our theory is not a limitation for applications.

The paper is organized as follows. The problem is stated in Section~\ref{stat} as a variational formulation of the PDE \eqref{e1}--\eqref{e1a} in appropriate function spaces depending on the choice of the boundary condition. A time discrete scheme with constant time step $\tau > 0$ is proposed in Section~\ref{disc}. Estimates independent of $\tau$ for the time discrete approximations are shown to be sufficient for passing to the limit in Section~\ref{proo} as $\tau \to 0$, and for proving that the limit is a solution to the original problem. Section~\ref{prop} is devoted to the proof of the bounded propagation speed of the moisture front when $p>3$, and to a discussion about the necessity of such a condition.


\section{Statement of the problem}\label{stat}

Given a parameter $\omega \in [0,1]$, we prescribe for $(x,t) \in \partial\Omega \times (0,T)$ the boundary condition
\be{e2}
	\omega \kappa(x,\theta) |\nabla u|^{p-2}\nabla u\cdot \mathbf{n}(x) + (1-\omega) u= 0 \ \mbox{ on } \partial \Omega.
\ee
The case $\omega = 1$ corresponds to the homogeneous Neumann boundary condition, $\omega = 0$ is the Dirichlet boundary condition, and $\omega \in (0,1)$ is the Robin boundary condition with boundary permeability $\gamma_\omega = (1-\omega)/\omega$. We define the reference space $X_\omega = W^{1,p}_\omega(\Omega)\cap L^\infty(\Omega)$, where $W^{1,p}_\omega(\Omega) = W^{1,p}(\Omega)$ for $\omega \in (0,1]$ and $W^{1,p}_\omega(\Omega)= W^{1,p}_0(\Omega)$ for $\omega = 0$, and state the problem in variational form
\begin{align}
	&\io \big(\theta_t \vrt + \kappa(x,\theta)|\nabla u|^{p-2}\nabla u\cdot\nabla \vrt\big)\dd x + \gamma_\omega \ipo u\vrt \dd s(x) = 0 \quad \forall \vrt \in X_\omega, \label{e4} \\
	&\theta = G[u] \ \ \ale \mbox{ in } \Omega\times(0,T), \label{e5} 
\end{align}
with the convention that $\gamma_0$ can be any real number.

The Preisach operator was originally introduced in \cite{prei}. For our purposes, it is convenient to use the equivalent variational setting from \cite{book}.

\begin{definition}\label{dpr}
Let $\lambda \in L^\infty(\Omega \times (0,\infty))$ be a given function with the following properties:
\begin{align}
	\label{ge6b}
	&\exists \Lambda>0 : \ \lambda(x,r) = 0\ \for r\ge \Lambda, \forall x \in \Omega,\\[2mm] \label{ge6}
	&\begin{aligned}
		\exists \bar{\lambda}>0 : \ |\lambda(x_1,r_1) - \lambda(x_2,r_2)| &\le \Big(\bar{\lambda}\,|x_1 - x_2| + |r_1 - r_2|\Big)\ \forall r_1, r_2 \in (0,\infty), \forall x_1,x_2 \in \Omega.
	\end{aligned}
\end{align}
For a given $r>0$, we call the {\em play operator with threshold $r$ and initial memory $\lambda$} the mapping which with a given function $u \in L^1(\Omega; W^{1,1}(0,T))$ associates the solution $\xi^r\in L^1(\Omega; W^{1,1}(0,T))$ of the variational inequality
\be{ge4a}
|u(x,t) - \xi^r(x,t)| \le r, \quad \xi^r_t(x,t)(u(x,t) - \xi^r(x,t) - z) \ge 0 \ \ale \ \forall z \in [-r,r],
\ee
with initial condition
\be{ge5}
\xi^r(x,0) = \lambda(x,r) \ \ale,
\ee
and we denote
\be{ge4}
\xi^r(x,t) = \play_r[\lambda,u](x,t).
\ee
Given a measurable function $\rho :\Omega\times(0,\infty)\times \real \to [0,\infty)$ and a bounded measurable function $\bar G :\Omega \to [0,1]$, the Preisach operator $G$ is defined as a mapping $G: L^2(\Omega; W^{1,1}(0,T))\to L^2(\Omega; W^{1,1}(0,T))$ by the formula
\be{ge3}
G[u](x,t) = \bar G(x) + \ \int_0^\infty\int_0^{\xi^r(x,t)} \rho(x,r,v)\dd v\dd r.
\ee
The Preisach operator is said to be {\em regular\/} if the density function $\rho$ of $G$ in \eqref{ge3} belongs to $L^\infty(\Omega\times (0,\infty )\times \real)$, and there exist constants $\rho_1,\bar{\rho}>0$ and a decreasing function $\rho_0: \real\to \real$ such that for all $U>0$, all $x,x_1,x_2 \in \Omega$, and a.\,e.\ $(r,v)\in (0,U) \times (-U,U)$ we have
\begin{align}
	&0 < \rho_0(U) < \rho(x,r,v) < \rho_1, \label{ge3a} \\[2mm]
	&|\rho(x_1,r,v) - \rho(x_2,r,v)| \le \bar{\rho}\,|x_1 - x_2|. \label{ge3b}
\end{align}
\end{definition}

In applications, the natural physical condition $\theta=G[u] \in [0,1]$ is satisfied for each input function $u$ if and only if the additional assumptions
\be{irho}
	\int_0^\infty \int_0^\infty \rho(x,r,v)\dd v\dd r \le 1-\bar G(x), \quad \int_0^\infty \int_0^\infty \rho(x,r,-v)\dd v\dd r \le \bar G(x),
\ee
hold for a.\,e.\ $x\in \Omega$. The function $\bar G(x)$ can be interpreted as residual moisture content at standard pressure. Note that the existence result in Theorem~\ref{t1} below holds independently of any specific choice of $\bar G$ and $\rho$ and, in particular, of \eqref{irho}.

Let us mention the following classical result (see \cite[Proposition~II.3.11]{book}).

\begin{proposition}\label{pc1}
	Let $G$ be a regular Preisach operator in the sense of Definition~\ref{dpr}. Then it can be extended to a Lipschitz continuous mapping $G: L^q(\Omega; C[0,T]) \to L^q(\Omega; C[0,T])$ for every $q \in [1,\infty)$.
\end{proposition}

\begin{figure}[htb]
	\begin{center}
		\includegraphics[width=.85\textwidth]{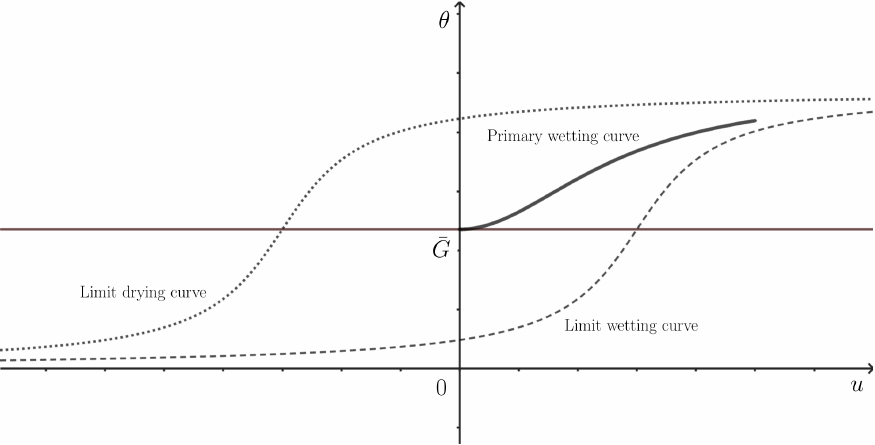}
		\caption{Primary wetting curve and limit wetting/drying curves in a typical Preisach diagram.}\label{fi0}
	\end{center}
\end{figure}

The Preisach operator is rate-independent. Hence, for an input function $u(x,t)$ which is monotone in a time interval $t\in (a(x),b(x))$, a regular Preisach operator $G$ can be represented by a superposition operator
\be{branch}
G[u](x,t) = G[u](x,a(x)) + B(x, u(x,t))
\ee
with an increasing function $u \mapsto B(x, u)$ called a {\em Preisach branch\/}. Indeed, the branches may be different at different points $x$ and different intervals $(a(x),b(x))$. The branches corresponding to increasing inputs are said to be {\em ascending\/} (the so-called wetting curves in the context of porous media), the branches corresponding to decreasing inputs are said to be {\em descending\/} (drying curves). In the ascending case, by \cite[Eq.~(2.25)]{ele} we have for $t \in (a(x),b(x))$ that
$$
\xi^r(x,t) = \max\{u(x,t) - r, \xi^r(x,a(x))\},
$$
hence the ascending branch $B_+(x,u)$ is given for $u \ge u(x,a(x))$ by the formula
\be{b+}
B_+(x,u) = \int_0^{r_0(x,u)}\int_{\xi^r(x,a(x))}^{u-r} \rho(x,r,v)\dd v\dd r,
\ee
where $r_0(x,u) = \min\{r>0: u-r \le \xi^r(x,a(x))\}$. Similarly, in the descending case it holds that
\be{b-}
B_-(x,u) = \int_0^{r_0(x,u)}\int^{\xi^r(x,a(x))}_{u+r} \rho(x,r,v)\dd v\dd r
\ee
for $u\le u(x,a(x))$ and $r_0(x,u) = \min\{r>0: u+r \ge \xi^r(x,a(x))\}$.
The wetting curve with initial memory $\lambda=0$ is called \textit{primary wetting curve}. All drying and wetting curves are bounded from below by the {\em limit wetting curve\/}, which is the theoretical wetting curve starting from the completely dry state, and from above by the {\em limit drying curve\/}, which is the theoretical drying curve starting from the completely wet state, see Figure~\ref{fi0}.

\begin{definition}\label{dpc}
	Let $U>0$ be given. A Preisach operator is said to be {\em uniformly counterclockwise convex on $[-U,U]$\/} if for all inputs $u$ such that $|u(x,t)|\le U$ a.\,e., all ascending branches are uniformly convex and all descending branches are uniformly concave.
	
	A regular Preisach operator $G$ is called {\em convexifiable\/} if for every $U>0$ there exist a uniformly counterclockwise convex Preisach operator $P$ on $[-U,U]$, positive constants $g_*(U),g^*(U),\bar{g}(U)$, and a twice continuously differentiable mapping $g:[-U,U] \to [-U,U]$ such that
\be{hg}
	g(0)=0, \quad 0<g_*(U) \le g'(u) \le g^*(U), \quad |g''(u)| \le \bar g(U)\ \ \forall u\in [-U,U],
\ee
	and $G = P\circ g$.
\end{definition}

A typical example of a uniformly counterclockwise convex operator is the so-called {\em Prandtl-Ishlinskii operator\/} characterized by positive density functions $\rho(x,r)$ independent of $v$, see \cite[Section~4.2]{book}. Operators of the form $P\circ g$ with a Prandtl-Ishlinskii operator $P$ and an increasing function $g$ are often used in control engineering because of their explicit inversion formulas, see \cite{al,viso,kk}. They are called the {\em generalized Prandtl-Ishlinskii operators\/} (GPI) and represent an important subclass of Preisach operators. Note also that for every Preisach operator $P$ and every Lipschitz continuous increasing function $g$, the superposition operator $G = P\circ g$ is also a Preisach operator, and there exists an explicit formula for its density, see \cite[Proposition~2.3]{error}. Another class of convexifiable Preisach operators is shown in \cite[Proposition~1.3]{colli}.

The technical hypotheses on the permeability function $\kappa$ can be stated as follows.

\begin{hypothesis}\label{hy2}
The permeability $\kappa:\Omega \times \real \to \real$ is a bounded Lipschitz continuous function, more precisely, there exist constants $\kappa_*, \kappa^*,\bar{\kappa}$ such that for all $\theta, \theta_1, \theta_2\in \real$ and all $x, x_1, x_2 \in \Omega$ we have
\be{hka}
0< \kappa_* \le \kappa(x,\theta) \le \kappa^*, \quad |\kappa(x_1,\theta_1) - \kappa(x_2,\theta_2)| \le \bar\kappa\big(|x_1 - x_2| + |\theta_1 - \theta_2|\big).
\ee
\end{hypothesis}

Note that even a local solution to Problem~\eqref{e1}--\eqref{e1a} may fail to exist if for example $\lambda(x,r) \equiv 0$ and $\!\dive \big(\kappa(x,\theta_0) |\nabla u_0|^{p-2}\nabla u_0\big)\ne 0$, and we need an initial memory compatibility condition which we state here following \cite{colli}. A more detailed discussion on this issue can be found in the introduction to \cite{colli}.

\begin{hypothesis}\label{hy1}
Let the initial memory $\lambda$ and the Preisach density function $\rho$ be as in Definition~\ref{dpr}, and let $\omega$ be the parameter in \eqref{e2}. The initial pressure $u_0$ belongs to $W^{2,\infty}(\Omega)$ and there exist a constant $L>0$ and a function $r_0 \in L^\infty(\Omega)$ such that, for $\Lambda>0$ as in \eqref{ge6b}, $\supess_{x\in \Omega}|u_0(x)| \le \Lambda$ and the following initial compatibility conditions hold:
	\begin{align} \label{c0}
		\lambda(x,0) &= u_0(x) \ \ale \textup{ in } \Omega,\\ \label{c0a}
		\theta_0(x) &= G[u](x,0) = \bar G(x) + \ \int_0^\infty\int_0^{\lambda(x,r)} \rho(x,r,v)\dd v\dd r \ \ale \textup{ in } \Omega,\\ \label{c1}
		&\hspace{-16mm}\frac1L \sqrt{\big|\dive \big(\kappa(x,\theta_0) |\nabla u_0|^{p-2}\nabla u_0\big)\big|} \le r_0(x) \le \Lambda \ \ale \textup{ in } \Omega, \\ \label{c2}
		-\frac{\partial}{\partial r} \lambda(x,r) &\in \sign\Big(\dive \big(\kappa(x,\theta_0) |\nabla u_0|^{p-2}\nabla u_0\big)\Big) \ \ale \textup{ in } \Omega \for r\in (0,r_0(x)), \\ \label{c2a}
		&\hspace{-17mm} \omega\kappa(x,\theta_0(x))|\nabla u_0|^{p-2}\nabla u_0\cdot \mathbf{n}(x) +(1-\omega) u_0(x) = 0 \ \ale \textup{ on } \partial\Omega.
	\end{align}
\end{hypothesis}

Unlike \cite{colli}, here we do not need to assume $\!\dive\!\big(\kappa(x,\theta_0(x))|\nabla u_0|^{p-2}\nabla u_0(x)\big) \in L^\infty(\Omega)$ since it follows from the fact that $u_0 \in W^{2,\infty}(\Omega)$ together with assumptions \eqref{ge6}, \eqref{ge3b}, and Hypothesis~\ref{hy2}. Our main existence result reads as follows.

\begin{theorem}\label{t1}
	Let $\omega \in [0,1]$ be given, let Hypotheses~\ref{hy2} and \ref{hy1} hold, and let $G$ be a convexifiable Preisach operator	in the sense of Definition~\ref{dpc}. Then there exists a solution $u \in L^\infty(\Omega\times (0,T))$ to Problem~\eqref{e4}--\eqref{e5}, \eqref{e3} such that $\nabla u \in L^p(\Omega\times (0,T);\real^N)$, and both $u_t$ and $\theta_t = G[u]_t$ belong to the Orlicz space $L^{\Phi_{log}}(\Omega\times (0,T))$ generated by the function $\Phi_{log}(v) = v\log(1+v)$. If moreover $\kappa = \kappa(x)$ is independent of $\theta$, then the solution is unique.
\end{theorem}

Basic properties of Orlicz spaces are summarized in \cite[Section~4]{perme}. For a more comprehensive discussion, we refer the interested reader to the monographs \cite{ada,rare}.

In the last part of the paper we address the problem of bounded speed of propagation, that is, whether solutions with compactly supported initial data preserve the property of compact support for all $t>0$. We show indeed that the wetting front propagates with asymptotically vanishing speed under the following assumptions.

	\begin{hypothesis}\label{hR}
		Let the initial memory $\lambda$ be as in Definition~\ref{dpr}. We assume that there exist constants $R_1>R_0>0$ and $\Lambda>0$ such that $\overline{B_{R_1}} \subset \Omega$ and
		\begin{align}\label{r0}
			u_0(x) = \lambda(x,r) = 0 \ & \for \ |x| \ge R_0, \ r>0,\\[1mm] \label{r0l}
			|u_0(x)|\le \Lambda, \ \lambda(x,r) \le (\Lambda - r)^+ \ &\for \ |x|<R_0, \ r>0.
		\end{align}
\end{hypothesis}

\begin{theorem}\label{t2}
	Let the assumptions of Theorem~\ref{t1} be satisfied, and additionally assume that Hypothesis~\ref{hR} holds. If $p > 3$, the permeability $\kappa>0$ is constant, and $\bar G$ as well as the function $\rho = \rho(r,v)$ in \eqref{ge3} are independent of $x$, then 
\begin{itemize}
\item[{\rm (i)}] There exists $t_1> 0$ such that the unique solution $u$ to Problem~\eqref{e4}--\eqref{e5}, \eqref{e3} from Theorem~\ref{t1} vanishes in $\Omega \setminus B_{R_1}$ for $t \in (0,t_1)$. Moreover, all solutions corresponding to different values of the parameter $\omega \in [0,1]$ coincide in $(0,t_1)$.
\item[{\rm (ii)}] This solution can be extended to $\real^N \times (0,\infty)$ with support contained in the ball $B_{R(t)}$ with radius $R(t) = R_0 + C_p t^{1/p}$ with a constant $C_p>0$ depending only on the data.
\end{itemize}
\end{theorem}

Roughly speaking, Theorem~\ref{t2}\,(i) states that the boundary condition on $\partial\Omega$ is not active until the moisture front reaches the boundary of $\Omega$.

An upper bound on the size of the support in terms of a power $1/p$ of time as in Theorem~\ref{t2}\,(ii) was obtained in \cite{diaz}. The proof relies on the comparison with a traveling wave solution of power-law type involving $t$ and $|x|$. Instead, to prove Theorem~\ref{t2} we construct a special family of dominant traveling wave solutions compatible with the hysteresis terms, see Section~\ref{prop}. Another argument is used in \cite[Chapter~VI]{dib}, where the comparison is made with the Barenblatt solution, whose support also grows as a power of time, but the power additionally depends on the dimension $N$.

Note that although Theorem~\ref{t2} is formulated in terms of the pressure $u$, it actually determines the maximal size of the wet area (i.\,e., the support of $\theta$). This is because no hysteresis effects can arise as long as $u$ remains equal to zero during the whole history. Therefore, we will continue to refer to the boundary of the support of $u$ as the moisture front.


\section{Time discretization}\label{disc}

Let $\omega \in [0,1]$ be given. We proceed as in \cite{colli}, choose a discretization parameter $n \in \nat$, define the time step $\tau = T/n$, and replace \eqref{e4}--\eqref{e5} with its time discrete system for the unknowns $\{u_i: i = 1,\dots,n\} \subset W^{1,p}_\omega (\Omega)$ of the form
\be{dis1}
\io \left(\frac1\tau(G[u]_i - G[u]_{i-1})\vrt + \kappa(x,G[u]_i)|\nabla u_i|^{p-2}\nabla u_i\cdot\nabla\vrt\right)\dd x + \gamma_\omega \ipo u_i \vrt \dd s(x) = 0
\ee
for every test function $\vrt \in X_\omega$. Here, the time-discrete Preisach operator $G[u]_i$ is defined for an input sequence $\{u_i : i\in\nat\cup\{0\}\}$ by a formula of the form \eqref{ge3}, namely,
\be{de3}
G[u]_i(x) = \bar G(x) + \ \int_0^{\infty}\int_0^{\xi^r_i(x)} \rho(x,r,v)\dd v\dd r,
\ee
where $\xi^r_i$ denotes the output of the time-discrete play operator
\be{de4}
\xi^r_i(x) = \play_r[\lambda,u]_i(x)
\ee
defined as the solution operator of the variational inequality
\be{de4a}
|u_i(x) - \xi^r_i(x)| \le r, \quad (\xi^r_i(x) - \xi^r_{i-1}(x))(u_i(x) - \xi^r_i(x) - z) \ge 0 \quad \forall i\in \nat \ \ \forall z \in [-r,r],
\ee
with a given initial condition
\be{de5}
\xi^r_0(x) = \lambda(x,r) \ \ale
\ee
similarly as in \eqref{ge4a}--\eqref{ge5}. Note that the discrete variational inequality \eqref{de4a} can be interpreted as weak formulation of \eqref{ge4a} for piecewise constant inputs in terms of the Kurzweil integral, and details can be found in \cite[Section~2]{ele}.

For each $i\in \{1,\dots,n\}$, there is no hysteresis in the passage from $u_{i-1}$ to $u_i$, so that \eqref{dis1} is a standard quasilinear elliptic equation. The existence of a solution $u_i \in W^{1,p}_\omega (\Omega)$ follows from a classical argument based on Fourier expansion into eigenfunctions of the Laplacian, Brouwer degree theory, and a homotopy argument similarly as in \cite[Section 3]{gravi}. An introduction to topological methods for solving nonlinear partial differential equations can be found in \cite[Chapter~V]{fuku}.


\subsection{Uniform upper bounds}\label{unif}

The first step is an $L^\infty$ bound on $u_i$. This is achieved by testing \eqref{dis1} by $\vrt =H_\ve(u_i - \Lambda)$, with $H_\ve$ being a Lipschitz regularization
\begin{equation}\label{Heav}
H_\ve (s) = \left\{
\begin{array}{ll}
	0 & \for s \le 0,\\
	\frac{s}{\ve} & \for s\in (0,\ve),\\
	1 & \for s \ge \ve,
\end{array}
\right.
\end{equation}
of the Heaviside function for some $\ve > 0$ and $\Lambda$ as in Hypothesis~\ref{hy1}. Arguing as in \cite[Section~2.1]{perme}, we obtain
\be{Lam}
|u_i(x)| \le \Lambda, \quad |\xi^r_i(x)| \le (\Lambda - r)^+
\ee
for a.\,e.\ $x \in \Omega$ and all $r\ge 0$ and $i\in \{0, 1,\dots,n\}$. In particular, even if we do not assume the a priori boundedness of $G$ as in \eqref{irho}, by assumption \eqref{ge3a} we obtain
\begin{equation}\label{GiL}
	 \quad |G[u]_i| \le C
\end{equation}
with a constant $C>0$ independent of $i$ and $\tau$.

We further test \eqref{dis1} by $\vrt = u_i$ and get for $\theta_i = G[u]_i$
\be{des2}
\frac1{\tau}\io (G[u]_i - G[u]_{i-1})u_i\dd x + \io \kappa(x,\theta_i)|\nabla u_i|^p\dd x + \gamma_\omega\ipo |u_i|^2\dd s(x) = 0
\ee
for all $i\in\{1,\dots,n\}$. We define the functions
\be{psi}
\psi(x,r,\xi) \coloneqq \int_0^\xi\rho(x,r,v)\dd v, \quad \Psi (x,r,\xi) \coloneqq \int_0^\xi v \rho(x,r,v)\dd v.
\ee
Choosing $z = 0$ in \eqref{de4a} and using the fact that the function $\psi$ is increasing, we obtain in both cases $\xi^r_i\ge \xi^r_{i-1}$ or $\xi^r_i\le \xi^r_{i-1}$ the inequalities
$$
(\psi(x,r,\xi^r_i)-\psi(x,r,\xi^r_{i-1})) u_i \ge (\psi(x,r,\xi^r_i)-\psi(x,r,\xi^r_{i-1}))\xi^r_i \ge \Psi(x,r,\xi^r_i)-\Psi(x,r,\xi^r_{i-1}).
$$
Identity \eqref{des2} together with Hypothesis~\ref{hy2} then yield
\be{energ}
\frac1\tau \io\int_0^\Lambda \big(\Psi(x,r,\xi^r_i)-\Psi(x,r,\xi^r_{i-1})\big)\dd r \dd x +\kappa_*\io |\nabla u_i|^p \dd x \le 0
\ee
for $i \in \{1,\dots,n\}$. Summing up over $i$ and exploiting the fact that, by the definition of $\xi^r_0$ in \eqref{de5} and the assumptions on $\rho$ and $\lambda$ in	Definition~\ref{dpr},
$$
\int_0^\Lambda \Psi(x,r,\xi^r_0(x)) \dd r = \int_0^\Lambda\int_0^{\lambda(x,r)} v\rho(x,r,v) \dd v\dd r \le \frac{\rho_1}{2} \int_0^\Lambda \lambda^2(x,r) \dd r \le C,
$$
and using \eqref{Lam} we get the estimate
\be{energy}
\max_{i=0, ..., n} \supess_{x \in \Omega} |u_i(x)| + \tau\sumiz \io |\nabla u_i|^p\dd x \le C
\ee
with a constant $C>0$ independent of $\tau$.


\subsection{Convexity estimate}\label{conv}

Recall that the operator $G$ is convexifiable in the sense of Definition~\ref{dpc}, that is, for every $U>0$ there exists a twice continuously differentiable mapping $g:[-U,U] \to [-U,U]$ such that $g(0) = 0$, $0 < g_* \le g'(u) \le g^* < \infty$, $|g''(U)| \le \bar{g}$, and $G$ is of the form
\be{ne0}
G = P \circ g,
\ee
where $P$ is a uniformly counterclockwise convex Preisach operator on $[-U,U]$. Let us fix $U$ from \eqref{Lam} and the corresponding function $g$.

We need to define a backward step $u_{-1}$ satisfying the strong formulation of \eqref{dis1} for $i=0$, that is,
\be{e7}
\frac1\tau(G[u]_0(x) - G[u]_{-1}(x)) = \dive \big(\kappa(x,G[u]_0(x)) |\nabla u_0|^{p-2}\nabla u_0\big) \ \mbox{ in } \Omega
\ee
with boundary condition \eqref{c2a}.
Repeating the argument of \cite[Proposition~3.3]{colli}, we use assumptions \eqref{ge3a} and \eqref{c2} to find for each $0<\tau <\rho_0(U)/2L^2$ functions $u_{-1}$ and $G[u]_{-1}$ satisfying \eqref{e7} as well as, thanks to \eqref{c1}, the estimate
\be{inim}
\frac1\tau |u_0(x) - u_{-1}(x)| \le C
\ee
with a constant $C>0$ independent of $\tau$ and $x$. The discrete equation \eqref{dis1} extended to $i=0$ has the form
\be{dp1}
\io \left(\frac1\tau(P[w]_i - P[w]_{i-1})\vrt + \kappa(x,\theta_i)|\nabla u_i|^{p-2}\nabla u_i\cdot\nabla\vrt\right)\dd x + \gamma_\omega \ipo u_i\vrt \dd s(x) = 0
\ee
with $w_i = g(u_i)$, $\theta_i = G[u]_i$ for $i\in \{0,1,\dots,n\}$ and for an arbitrary test function $\vrt \in X_\omega$. We proceed as in \cite{colli} and test the difference of \eqref{dp1} taken at discrete times $i+1$ and $i$
\begin{equation}\label{dp2} 
\begin{aligned}
 &\gamma_\omega \ipo (u_{i+1} -u_i)\vrt \dd s(x)	+ \io \bigg(\frac1\tau\big(P[w]_{i+1} - 2P[w]_i + P[w]_{i-1}\big)\vrt \\
	&\ + \Big(\kappa(x,\theta_{i+1})|\nabla u_{i+1}|^{p-2}\nabla u_{i+1}-\kappa(x,\theta_i)|\nabla u_i|^{p-2}\nabla u_i\Big) \cdot\nabla\vrt\bigg)\dd x = 0
\end{aligned}
\end{equation}
by $\vrt = f(w_{i+1} - w_i)$ with 
\be{dp2f}
f(w) \coloneqq \frac{w}{\tau + |w|}.
\ee
The boundary term gives a positive contribution which will be neglected. As for the hysteresis term on the left-hand side of \eqref{dp2}, following the same steps as in \cite[Section~2.2]{perme} we obtain
\be{dp4}
\frac1\tau \sumim (P[w]_{i+1} - 2P[w]_i + P[w]_{i-1})f(w_{i+1} {-} w_i)
\ge \frac{\beta}{4}\sumim |w_{i+1} {-} w_i|\log\left(1 + \frac{|w_{i+1} {-} w_i|}{\tau}\right) - C
\ee
with a constant $C>0$ independent of $\tau$. We thus get from \eqref{dp2}
\begin{align}\nonumber
&\sumim \io |w_{i+1} - w_i|\log\left(1 + \frac{|w_{i+1} {-} w_i|}{\tau}\right)\dd x\\ \label{dp5}
&\qquad + \sumim\io \Big(\kappa(x,\theta_{i+1})|\nabla u_{i+1}|^{p-2}\nabla u_{i+1}-\kappa(x,\theta_i) |\nabla u_{i}|^{p-2}\nabla u_i\Big)\cdot \nabla f(w_{i+1} - w_i)\dd x \le C
\end{align}
with a constant $C>0$ independent of $\tau$. Let us denote for simplicity $\kappa_i = \kappa(x,\theta_i)$. Then
\begin{align}\label{grad}
&\Big(\kappa_{i+1}|\nabla u_{i+1}|^{p-2}\nabla u_{i+1}-\kappa_i |\nabla u_{i}|^{p-2}\nabla u_i\Big)\cdot \nabla f(w_{i+1} {-} w_i)\\ \nonumber
&\ \, = f'(w_{i+1} {-} w_i)\Bigg(\Big(\big(\kappa_{i+1} {-} \kappa_{i}\big)|\nabla u_{i}|^{p-2}\nabla u_i + \kappa_{i+1}(|\nabla u_{i+1}|^{p-2}\nabla u_{i+1} {-}|\nabla u_{i}|^{p-2} \nabla u_i)\Big)\times\\ \nonumber
&\qquad \times \Big(\big(g'(u_{i+1}) {-} g'(u_{i})\big)\nabla u_i + g'(u_{i+1})(\nabla u_{i+1} {-}\nabla u_i)\Big)\Bigg)\\\nonumber
&\ \, = f'(w_{i+1} {-} w_i)\Bigg(\!g'(u_{i+1})\kappa_{i+1}(|\nabla u_{i+1}|^{p-2}\nabla u_{i+1} {-}|\nabla u_{i}|^{p} \nabla u_i)(\nabla u_{i+1} {-}\nabla u_i)\\ \nonumber
&\qquad + \big(g'(u_{i+1}) {-} g'(u_{i})\big)\!\big(\kappa_{i+1} {-} \kappa_{i}\big)\!|\nabla u_i|^{p} + g'(u_{i+1})\big(\kappa_{i+1} {-} \kappa_{i}\big) |\nabla u_{i}|^{p-2}\nabla u_i \cdot (\nabla u_{i+1} {-}\nabla u_i)\\ \nonumber
&\qquad + \kappa_{i+1}\big(g'(u_{i+1}) {-} g'(u_{i})\big)\nabla u_i \cdot (|\nabla u_{i+1}|^{p-2}\nabla u_{i+1} {-}|\nabla u_{i}|^{p-2} \nabla u_i)\Bigg).
\end{align}
The functions $\kappa$ and $g'$ are bounded and Lipschitz continuous, and
\begin{equation}\label{fprime}
	f'(w_{i+1} {-} w_i) = \frac{\tau}{(\tau + |w_{i+1} {-} w_i|)^2}.
\end{equation}
Moreover, since $\theta_i = P[w]_i$ admits a representation similar to \eqref{de3}, by Hypothesis~\ref{hy2}, estimate \eqref{GiL}, and the Lipschitz continuity of the time-discrete play implied by \eqref{de4a} we obtain
$$
|\kappa(x,\theta_{i+1}) - \kappa(x,\theta_{i})| \le \bar{\kappa}|\theta_{i+1} - \theta_i| \le C|w_{i+1} - w_i| \ \ale,
$$
whereas, from assumption \eqref{hg},
$$
|g'(u_{i+1}) {-} g'(u_{i})| \le \bar{g}(U)|u_{i+1} {-} u_i| \le \bar{g}(U)(g_*(U))^{-1}|w_{i+1} {-} w_i|.
$$
These considerations, as well as the elementary identity
\begin{equation}\label{elemen}
\begin{aligned}
	&\hspace{-6mm}(|\bsa|^{p-2} \bsa-|\bsb|^{p-2} \bsb) \cdot (\bsa - \bsb) = |\bsa|^p + |\bsb|^p - (|\bsa|^{p-2}+|\bsb|^{p-2})\, \bsa \cdot \bsb\\
	& = \frac12\left(|\bsa|^{p-2} + |\bsb|^{p-2}\right)|\bsa - \bsb|^2 + \frac12 \left(|\bsa|^{p-2} - |\bsb|^{p-2}\right)\left(|\bsa|^2 - |\bsb|^2\right),
\end{aligned}
\end{equation}
which holds for arbitrary vectors $\bsa, \bsb \in \real^N$, allow us to estimate each of the four terms inside the big brackets on the right-hand side of \eqref{grad} as follows:
\begin{align}\label{co1}
&g'(u_{i+1})\kappa_{i+1}(|\nabla u_{i+1}|^{p-2}\nabla u_{i+1} {-}|\nabla u_{i}|^{p-2} \nabla u_i)\cdot(\nabla u_{i+1} {-}\nabla u_i)\\ \nonumber
&\qquad \ge c (|\nabla u_{i+1}|^{p-2} + |\nabla u_{i}|^{p-2})|\nabla u_{i+1} {-}\nabla u_i|^2 \ge 0,\\[2mm]\label{co2}
& \big(g'(u_{i+1}) {-} g'(u_{i})\big)\!\big(\kappa_{i+1} {-} \kappa_{i}\big) |\nabla u_i|^{p} \le C |w_{i+1} {-} w_i|^2|\nabla u_i|^{p},\\[2mm] \label{co3}
& g'(u_{i+1})\big(\kappa_{i+1} {-} \kappa_{i}\big) |\nabla u_{i}|^{p-2}\nabla u_i \cdot (\nabla u_{i+1} {-}\nabla u_i)\\ \nonumber
& \qquad \le C |w_{i+1} {-} w_i|\big(|\nabla u_{i}|^{p-1} +|\nabla u_{i+1}|^{p-1} \big)||\nabla u_{i+1} {-}\nabla u_i|,\\[2mm] \label{co4}
& \kappa_{i+1}\big(g'(u_{i+1}) {-} g'(u_{i})\big)\nabla u_i \cdot (|\nabla u_{i+1}|^{p-2}\nabla u_{i+1} {-}|\nabla u_{i}|^{p-2} \nabla u_i) \\ \nonumber
& \qquad \le C|w_{i+1} {-} w_i|\big(|\nabla u_{i}|^{p-1} +|\nabla u_{i+1}|^{p-1} \big)|\nabla u_{i+1} {-}\nabla u_i|.
\end{align}
The terms on the right-hand side of \eqref{co3} and \eqref{co4} can be estimated using Young's inequality with $c$ as in \eqref{co1} and a possibly large constant $C$, namely
\begin{align*}
&|w_{i+1} {-} w_i|\big(|\nabla u_{i}|^{p-1} +|\nabla u_{i+1}|^{p-1} \big)|\nabla u_{i+1} {-}\nabla u_i|\\[1mm]
& \qquad \le C |w_{i+1} {-} w_i|^2\big(|\nabla u_{i}|^{p} +|\nabla u_{i+1}|^{p} \big) + \frac{c}{4}\big(|\nabla u_{i}|^{p-2} + |\nabla u_{i+1}|^{p-2}\big)|\nabla u_{i+1} {-}\nabla u_i|^2,
\end{align*}
where the second term can then be reabsorbed into \eqref{co1}. Hence, by \eqref{fprime}, for the left-hand side of \eqref{grad} we get the lower bound
\begin{align}\label{grad2}
&\Big(\kappa_{i+1}|\nabla u_{i+1}|^{p-2}\nabla u_{i+1}-\kappa_i |\nabla u_{i}|^{p-2}\nabla u_i\Big)\cdot \nabla f(w_{i+1} {-} w_i)\\ \nonumber
&\qquad \ge \frac{-\tau C |w_{i+1} {-} w_i|^2\big(|\nabla u_{i}|^{p} +|\nabla u_{i+1}|^{p} \big)}{(\tau + |w_{i+1} {-} w_i|)^2} \ge -\tau C \big(|\nabla u_{i}|^{p} +|\nabla u_{i+1}|^{p} \big).
\end{align}
As a consequence of \eqref{dp5}, \eqref{energy}, and \eqref{grad2}, and by assumption \eqref{hg} and the definition of $w_i$, we thus have the crucial estimate
\be{dp6}
\sumim \io |u_{i+1} - u_i|\log\left(1 + \frac{|u_{i+1} {-} u_i|}{\tau}\right)\dd x \le C\left(1+ \tau \sumiz \io|\nabla u_i|^p \dd x\right) \le C
\ee
with a constant $C>0$ independent of $\tau$.


\section{Proof of Theorem~\ref{t1}}\label{proo}

We construct the sequences $\{\hat u\on\}$ and $\{{\bar u}^{(n)}\}$ of approximations by piecewise linear and piecewise constant interpolations of the solutions to the time-discrete problem associated with the division $t_i=i\tau=iT/n$ for $i=0,1,\dots,n$ of the interval $[0,T]$. Namely, for $x \in \Omega$ and $t\in [t_{i-1},t_i)$, $i = 1,\dots,n$, we define 
\begin{align}\label{pl}
\hat u\on(x,t) &= u_{i-1}(x) + \frac{t-t_{i-1}}{\tau} (u_i(x) - u_{i-1}(x)),\\
\label{pco}
\bar u\on(x,t) &= u_{i}(x),
\end{align}
continuously extended to $t=T$. Similarly, we consider
\begin{align*}
\hat G\on(x,t) &= G[u]_{i-1}(x) + \frac{t-t_{i-1}}{\tau} (G[u]_i(x) - G[u]_{i-1}(x)),\\
\bar G\on(x,t) &= G[u]_i(x).
\end{align*}
Then \eqref{dis1} is of the form
\be{de0w}
\io \left(\hat G\on_t\vrt + \kappa(x,\bar G\on)\big|\nabla \bar u\on\big|^{p-2} \nabla \bar u\on \cdot\nabla\vrt\right)\dd x + \gamma_\omega \ipo \bar u\on\vrt \dd s(x) = 0
\ee
for every test function $\vrt \in X_\omega$. Our goal is to let $n \to \infty$ (or, equivalently, $\tau \to 0$) in \eqref{de0w}, and obtain in the limit a solution to \eqref{e4}--\eqref{e5}. The argument is based on the estimate which follows from \eqref{energy} and \eqref{dp6}, namely
\be{dp6a}
\supess_{(x,t) \in \Omega\times (0,T)}|\hat u\on(x,t)|+ \int_0^T\io |\hat u\on_t|\log\left(1 + \hat u\on_t\right)\dd x\dd t + \int_0^T\io|\nabla \hat u\on|^p \dd x \dd t \le K
\ee
with a constant $K>0$ independent on $n$. In \cite{perme}, we have proved the following results.

\begin{lemma}\label{lphi}
Let $\Phi_{log}$ be as in Theorem~\ref{t1}, and let $K>0$ be as in \eqref{dp6a}. Then there exists a function $\alpha:(0,1)\to (0,1)$ such that $\lim_{\tau \to 0} \alpha(\tau) = 0$ and for all $v \in (0,2K)$ and $\tau \in (0,1)$ we have
	$$
	\frac{\Phi_{log}(v)}{\tau\Phi_{log}\left(\frac{v}{\tau}\right)} = \frac{\log(1+v)}{\log\left(1+\frac{v}{\tau}\right)} \le \alpha(\tau).
	$$
\end{lemma} 

\begin{proposition}\label{pp1}
Let the sequence $u\on$ satisfy the condition \eqref{dp6a}. Then it is compact in the space $L^1(\Omega; C[0,T])$.
\end{proposition}

Proposition~\ref{pp1} and the Lipschitz continuity of $G$ in $L^1(\Omega;C[0,T])$ stated in Proposition~\ref{pc1} imply that, passing to a subsequence if necessary, $G[\hat u\on] \to G[u]$ in $L^1(\Omega;C[0,T])$. To estimate the difference $|\bar G\on(x,t) - G[\hat u\on]|$, let us define
\begin{align*}
\tilde \xi_r\on(x,t) &= \xi^r_{i-1}(x) + \frac{t-t_{i-1}}{\tau}\big(\xi^r_i(x) - \xi^r_{i-1}(x)\big) \quad \for (x,t)\in \Omega\times(t_{i-1}, t_i],\\[2mm]
\bar \xi_r\on(x,t) &= \xi^r_{i}(x) \quad \for (x,t)\in \Omega\times(t_{i-1}, t_i],\\[2mm]\hat\xi_r\on(x,t) &= \play_r[\hat u\on](x,t) \quad \for (x,t)\in \Omega\times[0,T],
\end{align*}
with $\xi^r_i$ defined in \eqref{de4}. By \eqref{ge4a} and \eqref{de4a} we have
\begin{align}\label{hxi}
\frac{\partial\hat \xi_r\on}{\partial t}\big(\hat u\on - \hat \xi_r\on - z\big) &\ge 0 \quad \ale,\\[2mm]\label{txi}
\frac{\partial\tilde \xi_r\on}{\partial t}\big(\bar u\on - \bar \xi_r\on - z\big) &\ge 0 \quad \ale,
\end{align}
for all $|z|\le r$. Putting $z= \bar u\on - \bar \xi_r\on$ in \eqref{hxi} and $z= \hat u\on - \hat \xi_r\on$ in \eqref{txi} and summing up the inequalities, we get
\be{htxi}
\frac12\frac{\partial}{\partial t}\left(\big|\hat \xi_r\on - \tilde \xi_r\on\big|^2\right) \le \frac{\partial}{\partial t}\big(\hat \xi_r\on - \tilde \xi_r\on\big)\left(\big(\hat u\on - \bar u\on\big) + \big(\bar \xi_r\on - \tilde \xi_r\on\big)\right).
\ee
Note that \eqref{de4a} implies
\be{tbxi}
|\bar \xi_r\on(x,t) - \tilde \xi_r\on(x,t)| \le \big|\xi^r_i(x) - \xi^r_{i-1}(x)\big| \le \big|u_i(x) - u_{i-1}(x)\big| \ \for\! \ale \ (x,t) \in \Omega\times (t_{i-1},t_i].
\ee
We have $|\frac{\partial}{\partial t}\hat \xi_r\on| \le |\frac{\partial}{\partial t}\hat u\on|$ a.\,e., and from \eqref{htxi}, \eqref{tbxi} we obtain
\be{htx2}
\frac{\partial}{\partial t}\left(\big(\hat \xi_r\on - \tilde \xi_r\on\big)^2\right)(x,t) \le \frac{8|u_i(x) - u_{i-1}(x)|^2}{\tau} \quad \for\! \ale \ (x,t) \in \Omega\times (t_{i-1},t_i].
\ee
Hence, by \eqref{htx2}, \eqref{Lam}, \eqref{dp6}, and Lemma~\ref{lphi},
\begin{align} \nonumber
\io \max_{t\in[0,T]}\big|\hat \xi_r\on - \tilde \xi_r\on\big|^2(x,t) \dd x &\le C\sumim \io|u_{i+1} - u_{i}|^2\dd x \\[2mm] \label{hti}
&\le C \sumim \io |u_{i+1} - u_{i}|\log(1 + |u_{i+1} - u_{i}|)\dd x \le C \alpha(\tau)
\end{align}
with a constant $C$ independent of $\tau$.
Additionally, \eqref{ge3a}, \eqref{Lam}, and \eqref{tbxi} yield
\begin{equation}\label{hGbG}
	|\hat G\on(x,t) - \bar G\on(x,t)| \le C\big|u_i(x) - u_{i-1}(x)\big|
\end{equation}
so that, combining \eqref{hti}--\eqref{hGbG}, we obtain
$$
\io \max_{t\in[0,T]}\big|G[\hat u\on] - \bar G\on\big|(x,t) \dd x + \io \max_{t\in[0,T]}\big|G[\hat u\on] - \hat G\on\big|(x,t) \dd x \le C \alpha^{1/2}(\tau).
$$
Hence, both $\bar G\on$ and $\hat G\on$ converge to $G[u]$ strongly in $L^1(\Omega\times (0,T))$. Furthermore, the sequence $\{\nabla \bar u\on\}$ is bounded in $L^p\big((\Omega \times (0,T); \real^N\big)$. Hence, $\nabla \bar u\on \to \nabla u$ weakly in $L^p\big((\Omega \times (0,T); \real^N\big)$. To pass to the limit in the elliptic term in \eqref{de0w}, we use a variant of the Minty trick. Notice first that $\vrt = \bar u\on - u +\delta y$ with $y \in L^p(0,T; X_\omega)$ and $\delta \in \real$ is an admissible test function in \eqref{de0w}, and we obtain
\be{mo1}
\begin{aligned}
&\io \left(\hat G\on_t(\bar u\on - u +\delta y) + \kappa(x,\bar G\on)\big|\nabla \bar u\on\big|^{p-2} \nabla \bar u\on \cdot\nabla(\bar u\on - u +\delta y)\right)\dd x\\
&\qquad + \gamma_\omega \ipo \bar u\on(\bar u\on - u +\delta y)\dd s(x) = 0.
\end{aligned}
\ee
By monotonicity (see also \eqref{elemen}) we have
\be{min2}
\hspace{-4mm}\kappa(x,\bar G\on)\left(\big|\nabla \bar u\on\big|^{p-2} \nabla \bar u\on - \big|\nabla (u - \delta y)\big|^{p-2} \nabla (u - \delta y)\right)\cdot\nabla(\bar u\on - u +\delta y) \ge 0
\ee
for all $\delta \in \real$ and $y \in L^p(0,T; X_\omega)$, hence
\begin{align}\nonumber
&\kappa(x, \bar G\on) \big|\nabla \bar u\on\big|^{p-2}\nabla \bar u\on\cdot\nabla(\bar u\on - u +\delta y)\\[2mm] \label{min1}
&\hspace{12mm} \ge \kappa(x, \bar G\on) \big|\nabla (u - \delta y)\big|^{p-2}\nabla (u - \delta y)\cdot\nabla(\bar u\on - u +\delta y),
\end{align}
and $\bar u\on(\bar u\on - u +\delta y) \ge (u-\delta y)(\bar u\on - u +\delta y)$. From \eqref{mo1} we obtain
\begin{align}\nonumber
&\gamma_\omega \ipo (u-\delta y)(\bar u\on - u +\delta y)\dd s(x) + \io \Big(\hat G\on_t(\bar u\on - u +\delta y)\\[2mm]\label{mo2}
&\hspace{15mm} + \kappa(x,\bar G\on) \big|\nabla (u - \delta y)\big|^{p-2}\nabla (u - \delta y)\cdot\nabla(\bar u\on - u +\delta y)\Big)\dd x\dd t \le 0.
\end{align}
The functions $\hat G\on_t$ are bounded the Orlicz space $L^{\Phi_{log}}(\Omega\times (0,T))$ generated by the function $\Phi_{log}$. It can be considered as the dual space to the small Orlicz space $L^{\Phi_{exp}}_\sharp(\Omega\times (0,T))$ generated by the function $\Phi_{exp}(v) = \expe^v - v - 1$. Hence, passing to a subsequence, we conclude that $\hat G\on_t$ converge weakly* to $G[u]_t$ in $L^{\Phi_{log}}(\Omega\times (0,T))$ (see \cite[Section~5]{perme} for more details).
Since $\hat u\on$ and $\bar u\on$ are uniformly bounded, by the Lebesgue Dominated Convergence Theorem they converge strongly in $L^{\Phi_{exp}}_\sharp(\Omega\times (0,T))$, and we get
\be{lim1}
\lim_{n\to\infty} \int_0^T\io \hat G\on_t (\bar u\on - u +\delta y) \dd x\dd t = \delta \int_0^T\io G[u]_t y \dd x\dd t.
\ee
To pass to the limit in the gradient term of \eqref{mo2}, we rewrite it as
$$
\begin{aligned}
	\kappa(x,\bar G\on) W \cdot\nabla(\bar u\on - u +\delta y) &= \kappa(x,G[u]) W \cdot\nabla(\bar u\on - u) + \delta\kappa(x,\bar G\on) W \cdot\nabla y \\
	&\quad + \big(\kappa(x,\bar G\on) - \kappa(x,G[u])\big) W \cdot\nabla(\bar u\on - u)
\end{aligned}
$$
with $W = \big|\nabla (u - \delta y)\big|^{p-2}\nabla (u - \delta y) \in L^{p'}(\Omega\times (0,T); \real^N)$, where $\frac{1}{p} + \frac{1}{p'} = 1$. The functions $\nabla \bar u\on$ weakly converge to $\nabla u$ in $L^{p}(\Omega\times (0,T); \real^N)$, hence
\be{lim2}
\lim_{n\to\infty} \int_0^T\io \kappa(x,G[u])\, W \cdot\nabla(\bar u\on - u) \dd x\dd t = 0.
\ee
We have, up to a subsequence, $\bar G\on \to G[u]$ strongly in every $L^q$ and pointwise almost everywhere, and it follows from the Lebesgue Dominated Convergence Theorem that
\be{lim3}
\lim_{n\to\infty} \delta\int_0^T\io \kappa(x,\bar G\on)\, W \cdot\nabla y \dd x\dd t = \delta\int_0^T\io \kappa(x,G[u])\, W \cdot\nabla y \dd x\dd t.
\ee
We estimate the remaining integral using H\"older's inequality as
\begin{align*}
&\int_0^T\io \big(\kappa(x,\bar G\on)-\kappa(x,G[u])\big) W \cdot\nabla(\bar u\on - u) \dd x\dd t \\[2mm]
&\hspace{15mm} \le C \left(\int_0^T\io|\kappa(x,\bar G\on)-\kappa(x,G[u])|^{p'}|W|^{p'} \dd x\dd t\right)^{1/{p'}},
\end{align*}
and using again the Lebesgue Dominated Convergence Theorem we conclude that
\be{lim4}
\lim_{n\to\infty} \int_0^T\io \big(\kappa(x,\bar G\on)-\kappa(x,G[u])\big) W \cdot\nabla(\bar u\on - u) \dd x\dd t = 0.
\ee
Summarizing the above computations and passing to the limit in the boundary term, we see that
$$
\begin{aligned}
&\delta \int_0^T\bigg(\io \left( G[u]_t y + \kappa(x, G[u])\big|\nabla (u - \delta y)\big|^{p-2} \nabla (u - \delta y) \cdot \nabla y\right)\dd x\\
& \qquad + \gamma_\omega \ipo (u-\delta y)y \dd s(x)\bigg)\dd t \le 0
\end{aligned}
$$
for every $\delta \in \real$. Dividing the above inequality by $|\delta|$ and letting $\delta$ tend to $0$, we complete the existence part of the proof of Theorem~\ref{t1}.

Assume now that $\kappa = \kappa(x)$ is independent of $\theta$, and let $u_1, u_2$ be two solutions of \eqref{e4}--\eqref{e5} with the same initial condition $u_0$ and the same initial memory $\lambda$. We use the Hilpert trick (see \cite{hilp}, see also \cite[Proposition~4.3]{colli}) and test the identity
\be{uni}
\io \Big((G[u_1] - G[u_2])_t \vrt + \kappa(x)\big(|\nabla u_1|^{p-2}\nabla u_1 - |\nabla u_2|^{p-2}\nabla u_2\big)\cdot\nabla \vrt\Big)\dd x = 0
\ee
by $\vrt = H_\ve(u_1 - u_2)$, where $H_\ve$ is the Lipschitz continuous approximation of the Heaviside function defined in \eqref{Heav}. We conclude that $u_1 \le u_2$ a.\,e., which implies uniqueness.


\section{Bounded speed of propagation}\label{prop}

This section is devoted to the proof of Theorem~\ref{t2}. The idea is to construct a class of explicit solutions $\usi$ to the PDE
\be{e0}
G[\usi]_t - \kappa \dive (|\nabla \usi|^{p-2} \nabla \usi) = 0,
\ee
with $p>2$ and for $(x,t) \in \real^N \times (0,\infty)$, which vanish on some half-space associated with a unit vector $\bsi \in \real^N$, and to show afterward that these solutions are dominant over the solution from Theorem~\ref{t1}.


\subsection{Traveling wave solutions}

Let $0<R_0<R_1$ be as in Hypothesis~\ref{hR}. We fix $R\in (R_0, R_1)$ and a unit vector $\bsi \in \real^N$, and define $\usi$ in the form of a traveling wave
\be{ae1}
\usi(x,t) = \uu_c\big(ct+ R - \bsi\cdot x\big)
\ee
with a suitable constant speed $c > 0$ of propagation, where $\uu_c : \real \to \real$ is a continuously differentiable function depending on $c$ such that $\uu_c(z) = 0$ for $z \le 0$, $\uu_c'(z) > 0$ for $z>0$. By definition, the function $\usi$ vanishes in the half-space $\bsi\cdot x \ge R + ct$. Moreover, for every $x \in \real^N$, the function $t \mapsto \usi(x,t)$ is nondecreasing, and $\usi(x,0) = \uu_c(R-\bsi\cdot x)$. We have
\be{ae2}
G[\usi](x,t) = \Gamma\big(\uu_c\big(ct+ R - \bsi\cdot x\big)\big),
\ee
where $\Gamma$ is the primary wetting curve associated with $G$, that is, according to \eqref{branch}--\eqref{b+} and the subsequent discussion,
\be{ilc}
\Gamma(u) = \bar G + \int_0^u\int_0^{u-r}\rho(r,v)\dd v\dd r =: \bar G + \Gamma_0(u),
\ee
with $\Gamma_0(0) = 0$. Furthermore, 
$$
\nabla \usi(x,t) = -\bsi \uu_c'\left(ct+ R - \bsi\cdot x\right), \quad |\nabla \usi|^{p-2} \nabla \usi = -\bsi \left(\uu_c'\left(ct+ R - \bsi\cdot x\right)\right)^{p-1}
$$
and
\be{prime}
\dive \big(|\nabla \usi|^{p-2} \nabla \usi\big) = \left(\left(\uu_c'\left(ct+ R - \bsi\cdot x\right)\right)^{p-1}\right)',
\ee
where the prime denotes the derivative with respect to $z$. Hence, in view of \eqref{ae2}, $\usi$ is a solution to \eqref{e0} if and only if
\be{ae3}
c\Gamma_0(\uu_c(z)) = \kappa\big(\uu_c'(z)\big)^{p-1}
\ee
for all $z\ge 0$, that is
\be{ae4}
\uu_c'(z) = \left(\frac{c}{\kappa}\Gamma_0(\uu_c(z))\right)^{1/(p-1)}, \qquad \uu_c(0) = 0.
\ee
If a solution to \eqref{ae4} exists, then it is unique. Since $\Gamma_0$ is nonnegative and increasing, we conclude that $\uu_c$ is increasing and convex, and $\lim_{z\to \infty} \uu_c(z) = \infty$.

A necessary and sufficient condition for the existence of a solution to \eqref{ae4} reads
\be{ns}
F(u^*) \coloneqq \int_0^{u^*} \Gamma_0(u)^{-1/(p-1)} \dd u < \infty \quad \forall u^*>0.
\ee
Then
\be{fc}
\uu_c(z) = F^{-1}(c^* z), \quad c^* \coloneqq \left(\frac{c}{\kappa}\right)^{1/(p-1)}.
\ee
From \eqref{ge3a} it follows that $\Gamma_0(u) \le \frac{\rho_1}{2} u^2$. On the other hand, we have for all $u^*>0$ and $u \in (0,u^*)$ that $\Gamma'_0(u) \ge u \rho_0(u^*)$, hence
\be{delta}
\Gamma_0(u) \ge \min\left\{\Gamma_0(u^*), \frac12 \rho_0(u^*) u^2\right\} \ \mbox{ for all } u>0.
\ee
By the monotonicity of $\Gamma_0$ and $\rho_0$, we can fix $u^*>0$ and $\rho_* > 0$ such that $\Gamma_0(u^*) \ge \rho_*$, $\rho_0(u^*) \ge 2\rho_*$, and obtain from the above computations that
\be{delt}
\rho_* \left(\frac{u}{1+u}\right)^2 \le \Gamma_0(u) \le \frac{\rho_1}{2} u^2 \quad \mbox{ for all } u>0.
\ee
The necessary and sufficient condition \eqref{ns} for the existence of traveling wave solutions thus holds if and only if $p>3$. Note that by virtue of \eqref{delt}, for $z$ close to $0$ we have $\uu_c(z) \approx z^{(p-1)/(p-3)}$, $\uu_c'(z) \approx z^{2/(p-3)}$, $(\uu_c'(z))^{p-1} \approx z^{2(p-1)/(p-3)}$, so that formula \eqref{prime} is meaningful.

Consider now the homogeneous Dirichlet problem on $B_{R_1}$
\begin{equation}\label{weakR}
	\int_{B_{R_1}} \Big(G[u]_t \phi + \kappa |\nabla u|^{p-2} \nabla u \cdot \nabla \phi\Big)\dd x = 0
\end{equation}
for every test function $\phi \in W^{1,p}_0(B_{R_1})\cap L^\infty(B_{R_1})$. The existence of a unique solution to \eqref{weakR} follows from Theorem~\ref{t1} for $\omega = 0$ and $\Omega = B_{R_1}$. We prove the following intermediate result.

\begin{proposition}\label{tt1}
Let Hypothesis~\ref{hR} be satisfied and let $p>3$. Then there exists $t_1 > 0$ such that the solution to Problem~\eqref{weakR}, \eqref{e3} has compact support in $B_{R_1}$ for all $t \in (0,t_1)$.
\end{proposition}

\begin{figure}[htb]
	\begin{center}
		\includegraphics[width=.4\textwidth]{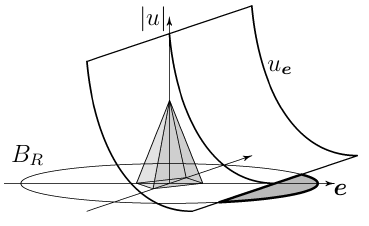}
		\caption{Dominant traveling wave solution.}\label{fi1}
	\end{center}
\end{figure}

\bpf{Proof}
The strategy is to compare the solution to Problem~\eqref{weakR}, \eqref{e3} with the traveling wave solution $\usi$ defined in \eqref{ae1}. We start by noting that, since $\uu_c$ is increasing and convex, we find $c>0$ such that
\be{alp}
\uu_c(R - R_0) \ge \Lambda.
\ee
This is indeed possible: by virtue of \eqref{fc}, condition \eqref{alp} can be equivalently written as
\begin{equation}\label{alp1}
\left(\frac{c}{\kappa}\right)^{1/(p-1)}(R - R_0) \ge F(\Lambda),
\end{equation}
which is certainly true if $c>0$ is sufficiently large. Then from \eqref{alp} and Hypothesis~\ref{hR} we obtain
$$
|u_0(x)|\le \uu_c(R - R_0) \le \uu_c(R - \bsi\cdot x) = \usi(x,0)
$$
for all $\bsi \in \partial B_1$ and a.\,e.\ $x \in B_{R_1}$. Let $t_1>0$ be such that $ct_1 + (R-R_0) < R_1$. For every test function $\phi \in W^{1,p}_0(B_{R_1})\cap L^\infty(B_{R_1})$, it follows from \eqref{e0} and \eqref{weakR} that
\be{test1}
\int_{B_{R_1}} \Big((G[u]- G[\usi])_t \phi + \kappa \big(|\nabla u|^{p-2} \nabla u- |\nabla \usi|^{p-2} \nabla \usi\big) \cdot \nabla \phi\Big)\dd x = 0.
\ee
Let $f:\real \to \real$ be a nondecreasing function. Then
$$
\big(|\nabla u|^{p-2} \nabla u- |\nabla \usi|^{p-2} \nabla \usi\big) \cdot \nabla f(u-\usi) \ge 0,
$$
see \eqref{elemen}. We then test \eqref{test1} by $H_\ve(u-\usi)$, where $H_\ve$ is the Lipschitz continuous approximation of the Heaviside function defined in \eqref{Heav}. This is indeed an admissible test function, since for $|x| = R_1$ we have by definition $u(x,t) = 0$ and $\usi(x,t) \ge 0$, hence $H_\ve(u-\usi)(x,t) = 0$. Passing to the limit as $\ve \to 0$, we obtain
\be{test}
\irn (G[u]- G[\usi])_t H(u-\usi)\dd x \le 0,
\ee
and from Hilpert's inequality (see \cite{hilp}) we conclude that $u(x,t) \le \usi(x,t)$ a.\,e.\ for each $\bsi \in \real^N$, $|\bsi| = 1$. By symmetry, we get $-u(x,t) \le \usi(x,t)$. For $t\in (0,t_1)$ and $\bsi = x/|x|$ we have in particular
\begin{equation}\label{RT}
	|u(x,t)| \le \uu_c\big(ct + R - |x|\big) \quad \mbox{and} \quad \uu_c\big(ct+ R - |x|\big) = 0 \ \for R_1 >|x|\ge ct_1 +R,
\end{equation}
see Figure~\ref{fi1}. Hence $u$ has compact support in $B_{R_1}$ for $t < t_1$, which completes the proof of Proposition~\ref{tt1}.

\epf


\subsection{Proof of Theorem~\ref{t2}}

Let $u_1$ be a solution to Problem~\eqref{weakR}, \eqref{e3} for $p>3$ and $t\in (0,t_1)$, and whose initial state satisfies Hypothesis~\ref{hR}. By Proposition~\ref{tt1}, it has compact support in $B_{R_1}$. If $u_1$ is extended by $0$ to $\Omega \setminus B_{R_1}$, then it necessarily coincides with the unique solution to Problem~\eqref{e4}--\eqref{e5}, \eqref{e3}, which is precisely statement~(i) of Theorem~\ref{t2}.

\begin{figure}[htb]
\begin{center}
\includegraphics[width=10cm]{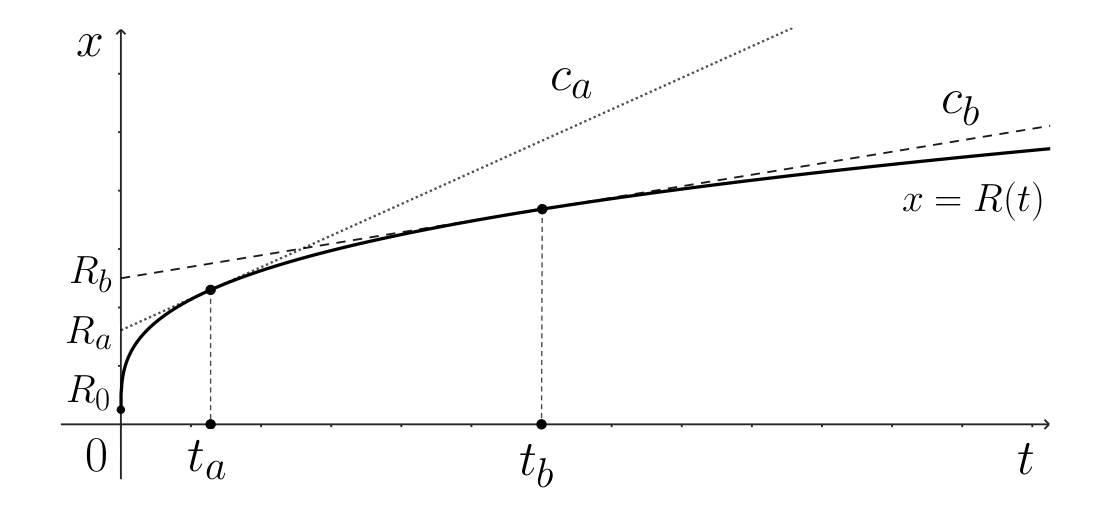}
\caption{Moving wetting front.}\label{fi2}
\end{center}
\end{figure}

To prove statement (ii), we repeat the argument of part (i), where we choose $R>R_0$ and $T>0$ arbitrarily large, and set $\Omega = B_{R_T}$ with $R_T \coloneqq cT + R$. According to formula \eqref{RT}, the corresponding solution $u_T(\cdot,t)$ has compact support in $B_{R_T}$ for $t\in (0,T)$ and coincides with $u_1$ for $t \in (0,t_1)$.
We can let $T\to \infty$ and conclude that $u_T$ can be extended to a solution $u_\infty$ on the whole space-time domain $\real^N \times (0,\infty)$.
According to \eqref{alp1}, $u_\infty$ is dominated by $\usi$ provided that
\be{fc1}
c^{1/(p-1)}(R - R_0) \ge \kappa^{1/(p-1)} F(\Lambda) \eqqcolon \bar\Lambda.
\ee
From \eqref{fc1} it follows that the speed $c>0$ of propagation can be arbitrarily small if $R$ is sufficiently large. An upper bound $x = R(t)$ for the wetting front is given by the envelope of the linear functions
\be{enve}
R_c(t) = R + ct, \quad R = R_0 + \bar\Lambda c^{-1/(p-1)},
\ee
defined by the condition that the straight lines $x = R_c(t)$ are tangent to the curve $x = R(t)$ at their contact points, see Figure~\ref{fi2}. This leads to the ODE
\be{ode}
R(t) - t \dot R(t) = R_0 + \bar\Lambda\big(\dot R(t)\big)^{-1/(p-1)},
\ee
which admits an explicit representation
$$
R(t) = R_0 + C_p\, t^{1/p}, \quad C_p = p\left(\frac{\bar\Lambda}{p-1}\right)^{(p-1)/p}.
$$
This completes the proof of Theorem~\ref{t2}.

\subsection{Doubly nonlinear equations and necessity of the condition $p>3$}\label{nece}

If $2<p\le 3$ the comparison technique developed in the previous two subsections cannot be applied. Obviously, this does not imply that for this range of exponents the speed of propagation is infinite. To fill this gap, in this subsection we gather some observations and draw some connections between our problem and some well-known nonlinear parabolic equations.

The double-sided estimate \eqref{delt} indicates that Problem~\eqref{e1}--\eqref{e1a} can be compared to the equation
\be{dne1}
\big(|u|u\big)_t - \kappa\dive\big(|\nabla u|^{p-2}\nabla u\big) = 0,
\ee
which is an instance of the doubly nonlinear parabolic equation
\be{dne2}
\big(|u|^{m-1}u\big)_t - \kappa\dive\big(|\nabla u|^{p-2}\nabla u\big) = 0
\ee
when $m=2$. The case $p-1>m$ is commonly known as the \textit{slow diffusion equation}, while the case $p-1<m$ is named the \textit{fast diffusion equation}, see \cite{iva,boge}. The difference between them becomes apparent in the fact that slow diffusion equations allow solutions with compact support and perturbations propagate at finite speed, while in fast diffusion equations perturbations propagate at infinite speed, prohibiting solutions with compact support. When $m=2$ as in our case, the critical exponent is exactly $p=3$. On the one hand, this confirms that when $p>3$ the speed of propagation has to be finite, which is in agreement with our findings; on the other hand, when $p<3$ (or $p=3$) we cannot expect solutions with compactly supported initial data to be compactly supported for any positive time, so the condition $p>3$ is also necessary.


\end{document}